\documentclass[a4paper,12pt]{amsart}
\usepackage{amsmath}
\usepackage{amsfonts}
\usepackage{amscd}
\usepackage{amssymb}

\begin{document}


\font\msbm=msbm10
\newcommand{\B}[1]{\hbox{\msbm #1}}
\newcommand{\C}[1]{\hbox{\Cal #1}}
\newcommand{\F}[1]{\frak #1}

\newcommand{\MS}{{\medskip}}
\newcommand{\SK}{{\smallskip}}
\newcommand{\BS}{{\bigskip}}
\newcommand{\NI}{{\noindent}}

\newtheorem{theorem}{Theorem}[section]
\newtheorem{thm}[theorem]{Theorem}
\newtheorem{cor}[theorem]{Corollary}
\newtheorem{defin}[theorem]{Definition}
\newtheorem{rem}[theorem]{Remark}
\newtheorem{lemma}[theorem]{Lemma}
\newtheorem{prop}[theorem]{Proposition}
\newtheorem{ex}[theorem]{Example}
\newtheorem{que}[theorem]{Question}
\newtheorem{con}[theorem]{Conjecture}
\numberwithin{equation}{section}
\newtheorem{nta}[theorem]{Notation}

\title{\bf The Lalonde-McDuff conjecture for nilmanifolds}
\author{Zofia St\c{e}pie\'{n}}
\date{August 8, 2006.}
\maketitle

\begin{abstract}
\NI We prove that any Hamiltonian bundle whose fiber is a
nilmanifold c-splits.
\smallskip

\noindent {\bf Keywords:}  Hamiltonian fiber bundle, c-splitting,
minimal Sullivan model.

\noindent {\bf AMS classification(2000)}: 55P62, 57R17.

\end{abstract}

\section{Introduction}
Let $(M,\omega)$ be a closed symplectic manifold. A fibre
bundle $(M,\omega)\rightarrow P\rightarrow B$ is called
{\em Hamiltonian} if its structural group can be reduced
to the group of Hamiltonian diffeomorphisms.
The following conjecture was posed by Lalonde and McDuff in
\cite{LM}

\begin{con}
Every Hamiltonian fiber bundle c-splits. This
means that there is an additive isomorphism
$$H^{*}(P)\cong H^{*}(B)\otimes H^{*}(M)$$
\end{con}
where $H^{*}$ denotes cohomology with real coefficients.

\SK

The c-splitting conjecture holds in many cases, see \cite{H},
\cite{K}, \cite{LM}, but the general case is still not resolved.
In particular, an argument of Blanchard, \cite{B} shows that it
holds when the cohomology of the fiber satisfies the hard
Lefschetz condition. The purpose of this paper is to prove the
c-splitting conjecture for fibers which are symplectic
nilmanifolds.  A {\em nilmanifold} is a compact homogeneous space
of the form $N/\Gamma$, where $N$ is a simply connected nilpotent
Lie group and $\Gamma$ is a discrete co-compact subgroup (i.e. a
lattice). Let us recall that these manifolds do not satisfy the
hard Lefschetz condition except for tori (\cite{BG}, \cite{M1},
\cite{TO}).

In this
article we prove the following theorem.

\begin{thm}
Let $(N/\Gamma,\omega)$ be a symplectic nilmanifold and
$(N/\Gamma,\omega)\hookrightarrow P\rightarrow B$ be a Hamiltonian
bundle over a simply connected  CW-complex. Then $H^{*}(P)\cong
H^{*}(B)\otimes H^{*}(N/\Gamma)$ as algebras.
\end{thm}

This theorem implies the following result.

\begin{thm}
Let $(N/\Gamma,\omega)$ be a symplectic nilmanifold. Then any
Hamiltonian bundle $(N/\Gamma,\omega)\hookrightarrow P\rightarrow
B$ c-splits.
\end{thm}

\begin{proof} Since $BHam(N/\Gamma)$ is simply connected then the universal Hamiltonian bundle with fiber $N/\Gamma$
$$ N/\Gamma \rightarrow
EHam(N/\Gamma)\times_{Ham(N/\Gamma)}N/\Gamma\rightarrow
BHam(N/\Gamma)$$ is c-split by Theorem 1.2. It follows from
\cite[Lemma 4.1(i)]{LM} that any Hamiltonian bundle with fiber
$N/\Gamma$ c-splits.
\end{proof}

We will prove  Theorem 1.2 in the next section. The first step is
to establish the result for $B=S^{2}.$ In \cite{M2},  McDuff
proved that every Hamiltonian bundle over $S^{2}$ c-splits. The
proof is difficult and requires hard analytic tools. In the
special case of nilmanifold fiber our approach gives a different
and easy proof. The main tool used in this paper is the Sullivan
model of a fibration. Also, we use the following characterization
of Hamiltonian bundles:

\begin{thm}\cite[Lemma 2.3]{LM}
If $\pi_{1}(B)=0$ then a symplectic bundle $M\hookrightarrow
P\rightarrow B$ is Hamiltonian if and only if the class
$[\omega]\in H^{2}(M)$ extends to $a\in H^{2}(P)$.
\end{thm}

\bigskip
\section{Proof of Theorem 1.2}

\SK

\subsection{Some facts from rational homotopy theory}
Basic reference for  rational homotopy theory is \cite{FHT}. We
also keep terminology and notation close to \cite{FHT}. In this
paper all algebras and cohomologies are considered over the field
of real numbers.

A free commutative cochain algebra $(\Lambda V,d)$ such that
\begin{enumerate}
    \item $V$ is a graded vector space and $V=\{V^{p}\}_{p\geq1}$,
    \item $V=\bigcup^{\infty}_{k=0}V(k)$, where
    $V(0) \subset V(1) \subset \cdots$ is an increasing sequence of
    graded subspaces such that
    $$d=0 \,\, \text{in} \,\, V(0)\,\,\,\, \text{and} \,\,\,\,
    d:V(k)\rightarrow\Lambda V(k-1),\,\,k\geq1,$$
    \item $\text{Im}\, d \subset\Lambda^{+}V\cdot\Lambda^{+}V$,
    where $\Lambda^{+}V=\oplus^{\infty}_{q=1}\Lambda^{q}V$ and
$\Lambda^{q}V$ is the linear span of
    the elements $v_{1} \wedge \cdots \wedge v_{q}$, $v_{i} \in
    V,$
    \item there is a quasi-isomorphism
    $$m:(\Lambda V,d)\rightarrow A_{PL}(X).$$
\end{enumerate}
is called the {\em minimal Sullivan model for a path connected
topological space $X$.} Here $A_{PL}(X)$ denotes the cochain
algebra of polynomial differential forms on $X.$

If condition (3) is not satisfied then $(\Lambda V,d)$ is called
the {\em Sullivan model for $X$}. For any path connected
topological space $X$ there exists a minimal Sullivan model
$(\Lambda V,d)$ and this is uniquely determined up to isomorphism.
If $X$ is a finite CW complex then each $V^{p}$ is finite
dimensional.

\SK

Let $f:X\rightarrow Y$ be a continuous map between path connected
topological spaces. A choice of $m:(\Lambda V_{Y},d)\rightarrow
A_{PL}(Y)$ and $m':(\Lambda V_{X},d')\rightarrow A_{PL}(X)$
determines a unique homotopy class of morphisms between Sullivan
models of $Y$ and $X$ respectively. A morphism $f^{*}:(\Lambda
V_{Y},d)\rightarrow (\Lambda V_{X},d')$ such that $m'f^{*}$ and
$A_{PL}(f)m$ are homotopic is called a {\em Sullivan
representative of f}. Here $A_{PL}(f):A_{PL}(Y)\rightarrow
A_{PL}(X)$ is a morphism induced by $f.$

\SK

The following Theorem can be found in \cite{AP}. See also
 Proposition 15.5 in \cite {FHT}.

\begin{thm}{(Grivel-Halperin-Thomas)} Let $\pi:E \rightarrow B$ be a
Serre fibration of path connected spaces and $F=\pi^{-1}(b)$ be
the fiber over the base-point $b.$ Suppose that:
\begin{enumerate}
    \item $F$ is path connected,
    \item $\pi_{1}(B)$ acts nilpotently on $H^{k}(F)$ for all
    $k\geq1,$
    \item either $B$ or $F$ is of finite type.
\end{enumerate}
Then there exist a commutative diagram of cochain algebra
morphisms
$$
\begin{CD}
A_{PL}(B) @>>> A_{PL}(P) @>>> A_{PL}(F) \\
@Am_{B}AA             @AmAA         @Am_{F}AA \\
(\Lambda V_{B},d)  @>>> (\Lambda V_{B}\otimes \Lambda V_{F} , D)
@>>>
(\Lambda V_{F}, \overline{d})  \\
\end{CD}
$$
in which $(\Lambda V_{B},d)$ is minimal Sullivan model for $B$,
$(\Lambda V_{F}, \overline{d})$ is minimal Sullivan model for $F$,
$m$ is a quasi-isomorphism and

$$Db=db,\,\, b\in
\Lambda V_{B}$$
$$Dv-\overline{d}v \in \Lambda^{+} V_{B}\otimes \Lambda V_{F},\,\,v\in V_{F}.$$

\end{thm}

\subsection{A minimal Sullivan model for $S^{2}$.}
The fundamental class $[S^{2}]\in H_{2}(S^{2};\mathbb{Z})$
determines a class $\omega\in H^{2}(A_{PL}(S^{2}))$ such that
$\langle\omega, [S^{2}]\rangle=1.$ Here $\langle,\rangle$ denotes
the paring between cohomology and homology. Let $\Phi$ be a
cocycle representing the cohomology classes $\omega$ and
$\Phi^{2}=d\Psi.$ Then the a minimal Sullivan model for $S^{2}$ is
given by
$$m:(\Lambda(a,b),d)\rightarrow A_{PL}(S^{2})$$
where $deg(a)=2,$ $deg(b)=3$, $ma=\Phi$, $mb=\Psi$

\subsection{A minimal Sullivan model for a nilmanifold.}

Let $N/\Gamma$ be a nilmanifold of dimension $n$. In this case,
the differential graded algebra $A_{DR}^{N}(N)$ of right-invariant
forms on $N$ is a Sullivan minimal model for $N/\Gamma$. In other
words, there is a dual basis $\{x_{1},x_{2},...,x_{n}\}$ to the
basis $\{X_{1},X_{2},...,X_{n}\}$ of ${\F n}$ such that

$$(\Lambda {\F n}^{*},\overline{d})=(\Lambda(x_{1},x_{2},...,x_{n}),\overline{d})$$
with $deg(x_{k})=1$, $k\in\{1,2,...,n\}$ is a minimal model of
$N/\Gamma$. Here ${\F n}$ denotes the Lie algebra of the Lie group
$N$.

\SK

 We call a form
$\alpha$ on $N/\Gamma$ homogeneous if the pullback of $\alpha$ to
$N$ is right-invariant. Thus, any cohomology class $[\alpha]\in
H^{*}(G/\Gamma)$ is represented by a homogeneous form
$\alpha_{h}.$
  Recall (\cite{M1}, \cite{IRTU}) that the cohomology  class of the symplectic structure is
represented by a homogeneous symplectic form. Therefore, it is
represented by a degree 2 element of the minimal model which we
can write as $\omega=\Sigma a_{ij}x_{i}x_{j}$ and the following
conditions are satisfied:
\begin{enumerate}
    \item $\overline{d}\omega=0,$
    \item for every $X \in {\F n}$\\$\omega(X,-)=0\,\,\Rightarrow\,\, X=0$
\end{enumerate}

\SK

\subsection{Proof of Theorem 1.2} Let $N/\Gamma$ be a symplectic
nilmanifold of dimension $2n$ and
$$ N/\Gamma\hookrightarrow P\rightarrow B$$ be a Hamiltonian
bundle over a simply connected CW-complex. The assumptions of
Theorem 2.1 are satisfied since $N/\Gamma$ is of finite type. Thus
we obtain  a commutative diagram

$$
\begin{CD}
A_{PL}(B) @>>> A_{PL}(P) @>>> A_{PL}(N/\Gamma) \\
@Am_{B}A \simeq A             @Am_{P}A \simeq A         @Am_{N/\Gamma}A \simeq A \\
(\Lambda V_{B},d)  @>>> (\Lambda V_{B}\otimes \Lambda {\F n}^{*} ,
D) @>>>
(\Lambda {\F n}^{*}, \overline{d})  \\
\end{CD}
$$

\noindent in which $(\Lambda V_{B},d)$ is a minimal Sullivan model
for $B$, $(\Lambda {\F n}^{*},\overline{d})$ is a minimal Sullivan
model for $N/\Gamma$ and

$$Db=db,\,\, b\in
\Lambda V_{B}$$
$$Dv-\overline{d}v \in \Lambda^{+} V_{B}\otimes \Lambda {\F
n}^{*},\,\,v\in {\F n}^{*}.$$

\noindent
We will show that $Dv=\overline{d}v,\,\,v\in {\F n}^{*}$. This
implies that $H^{*}(P)= H^{*}(B)\otimes H^{*}(N/\Gamma)$

\SK

\noindent
{\bf Case 1.} $B=S^{2}$. Let
$((\Lambda(a,b)\otimes\Lambda(x_{1},x_{2},...,x_{2n}),D)$ with
$deg(x_{k})=1$, $k\in\{1,2,...,2n\}$ be the  Sullivan model for
$P$.
By degree reasons the differential $D$ is given by

$$D(x_{k})=\alpha_{k}a+\overline{d}x_{k},\,\,\text{where}\,\, \alpha_{k}\in \mathbb{R}$$

\SK

 Let $\omega=\sum a_{ij}x_{i}x_{j}$ represent the cohomology class of the symplectic form on
 $N/\Gamma$ in the minimal model
 $\Lambda(x_{1},x_{2},...,x_{2n}).$
 Since the bundle is Hamiltonian we have  $D(\omega)=0$ (Theorem
 1.4).
 \SK

Thus
\begin{eqnarray*}
0&=&D(\omega)\\
&=&D(\sum a_{ij}x_{i}x_{j})\\
&=&\sum a_{ij}(Dx_{i}x_{j}-x_{i}Dx_{j})\\
&=&\sum a_{ij}((\alpha_{i}a+\overline{d}x_{i})x_{j}-
                x_{i}(\alpha_{j}a+\overline{d}x_{j}))\\
&=&\sum a_{ij}(\alpha_{i}ax_{j}-\alpha_{j}x_{i}a)+\overline{d}\omega
\end{eqnarray*}

\noindent
and we obtain

$$\sum a_{ij}(\alpha_{i}x_{j}-\alpha_{j}x_{i})=0.\eqno (1)$$

\SK
Take $X\in {\F n}$ such that $x_{k}(X)=\alpha_{k}.$ Note that

$$\omega(X,-)=\sum a_{ij}(\alpha_{i}x_{j}-\alpha_{j}x_{i}).$$

\SK

\noindent It follows from (1) that

$$X=0,\, \text{and thus}\,  \alpha_{k}=0,\, \text{for
all}\, k.$$

\SK
\noindent
Thus $D(x_{k})=\overline{d}x_{k}$ and $H^{*}(P)=
H^{*}(S^{2})\otimes H^{*}(N/\Gamma)$.

\SK \noindent {\bf Case 2.} $B$ is a simply connected finite CW
complex. First, we will choose a minimal Sullivan model for $B.$
Let $v_{1},...,v_{m}$ be a basis of $H_{2}(B,\mathbb{R}).$ Since
$B$ is simply connected the Hurewicz Theorem asserts that the
Hurewicz map is an isomorphism $\pi_{2}(B)\rightarrow
H_{2}(B,\mathbb{Z}).$ This extends to an isomorphism
$\pi_{2}(B)\otimes_{\mathbb{Z}}\mathbb{R}\rightarrow
H_{2}(B,\mathbb{R})$. Thus  $v_{i}$ can be represented by maps
$f_{i}:S^{2}\rightarrow B.$ Let $a_{1},...,a_{m}$ be the basis for
$H^{2}(B,R)$ defined by $\langle
a_{j},H_{*}f_{i}([S^{2}])\rangle=\delta_{ij}$, where
$H_{*}f_{i}:H_{*}(S^{2})\rightarrow H_{*}(B)$ and
$\langle,\rangle$ denotes the paring between cohomology and
homology. Let $\Phi_{j}\in A_{PL}(B)$ be the cocycle representing
the cohomology class $a_{j}$. Since $V^{2}_{B}\simeq H^{2}$ then
we can choose $V_{B}$ and a quasi-isomorphism

$$m:(\Lambda V_{B},\widetilde{d})\rightarrow A_{PL}(B),$$
such that $ma_{j}=\Phi_{j},$ $j\in\{1,\ldots, m \}.$ Hence
$f_{i}^{*}(a_{j})=\delta_{ij} a$ for any $f_{i}^{*}.$  Here
$f_{i}^{*}:(\Lambda
V_{B},\widetilde{d})\rightarrow(\Lambda(a,b),d)$ is a Sullivan
representative of $f_{i}.$

Now, let

$$(\Lambda(V_{B})\otimes\Lambda(x_{1},x_{2},...,x_{2n}),\widetilde{D})$$
with $deg(x_{k})=1$, $k\in\{1,2,...,2n\}$ be the  Sullivan model
for $P$. Again, from degree reasons the differential $\widetilde{D}$ is
given by

$$\widetilde{D}(x_{k})=\sum \alpha_{kj}a_{j}+\overline{d}x_{k},$$

\noindent where $deg(a_{j})=2,\,\,\alpha_{kj}\in \mathbb{R}$.

\SK

 Now, consider pullback of the bundle to $S^{2}$ by the map
$f_{i}$ for each $i\in\{1,2,...,m\}.$ Then we can write the
commutative diagram

$$
\begin{CD}
(\Lambda(x_{1},x_{2},...,x_{2n}),\overline{d})
@=(\Lambda(x_{1},x_{2},...,x_{2n}),\overline{d})\\
             @AAA                  @AAA\\
(\Lambda(a,b)\otimes\Lambda(x_{1},x_{2},...,x_{2n}),D)           @<\widehat{f}_{i}^{*}<<   (\Lambda(V_{B})\otimes\Lambda(x_{1},x_{2},...,x_{2n}),\widetilde{D})) \\
             @AAA              @AAA   \\
             (\Lambda(a,b),d)          @<f_{i}^{*}<<     (\Lambda V_{B},\widetilde{d}) \\
\end{CD}
$$

\SK
It follows from Case 1 and the commutativity of the diagram

\begin{eqnarray*}
\widehat{f}_{i}^{*}(\widetilde{D}(x_{1}))&=&
\widehat{f}_{i}^{*}(\sum \alpha_{1j}a_{j}+\overline{d}x_{1})\\
&=&\sum \alpha_{1j}f_{i}^{*}(a_{j})+\overline{d}x_{1}\\
&=&\sum \alpha_{1j}(\delta_{ij}a)+\overline{d}x_{1}\\
&=&\alpha_{1i}a+\overline{d}x_{1}\\
&=&D(\widehat{f}_{i}^{*}(x_{1}))\\
&=&D(x_{1})=\overline{d}x_{1}
\end{eqnarray*}

\SK
\noindent
Hence $\alpha_{1i}=0,$ for all $i\in\{1,2,...,m\}$ The above
argument applied to $x_{2},x_{3},...,x_{2n}$ gives
$\widetilde{D}(v)=\overline{d}v,\,\,v\in {\F n}^{*}.$

\SK
\noindent
{\bf Case 3.}
$B$ is a simply connected CW-complex. A Sullivan model for $P$ is of the form
$(\Lambda V_{B}\otimes \Lambda {\F n}^{*},D)$,
where $\Lambda V_{B}$ is a minimal Sullivan model
for $B$. Assume that the condition $Dv=\overline{d}v,\,\,v\in {\F
n}^{*}$ is not satisfied. Then we can choose a finite subcomplex
$B_{1}\subset B$ such that the restriction of the bundle $P$ over
$B_{1}$ has a Sullivan model of the form $(\Lambda
V_{B_{1}}\otimes \Lambda {\F n}^{*},D)$, $V_{B_{1}} \subset V_{B}$
and the condition $Dv=\overline{d}v,\,\,v\in {\F n}^{*}$ is not
satisfied. However, this contradicts Case 2.
\qed

\bigskip

\NI {\bf Acknowledgements.} The author is grateful to Boguslaw
Hajduk and Jarek K\c{e}dra for valuable advice and continuous
support.

\bibliographystyle{amsalpha}

\medskip
\noindent Zofia St\c{e}pie\'{n}

\noindent Institute of Mathematics

\noindent
Szczecin University of
Technology

\noindent Al. Piast\'ow 17,

\noindent
70-310 Szczecin, Poland

\noindent \tt stepien@ps.pl

\end{document}